\documentclass{article}

\usepackage{arxiv}

\usepackage[utf8]{inputenc} 
\usepackage[T1]{fontenc}    
\usepackage{nameref,hyperref} 
\usepackage[dvipsnames,svgnames]{xcolor}
\hypersetup{
  colorlinks=true,
  citecolor=Mulberry,
  urlcolor=Bittersweet
}
\usepackage{url}            
\usepackage{booktabs}       
\usepackage{amsfonts}       
\usepackage{nicefrac}       
\usepackage{microtype}      
\usepackage{lipsum}         
\usepackage{graphicx}
\usepackage[square,sort,comma,numbers]{natbib}
\usepackage{doi}

\usepackage[labelfont=bf]{caption}
\usepackage{subcaption}
\usepackage{amsmath, amssymb}
\usepackage{cleveref}       
\usepackage{float}

\DeclareCaptionLabelFormat{custlettering}{{\Large \MakeUppercase{#2}}} 
\newenvironment{adjustwidth}[2]{%
}{%
    \ignorespacesafterend%
}

\title{Evaluation of data driven low-rank matrix factorization for accelerated solutions of the Vlasov equation}

\date{}

\newif\ifuniqueAffiliation
\uniqueAffiliationtrue

\ifuniqueAffiliation 
\author{ \href{https://orcid.org/0009-0009-1054-0995}{\includegraphics[scale=0.06]{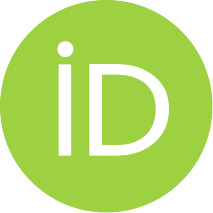}\hspace{1mm}Bhavana Jonnalagadda}\thanks{Code available at \href{https://github.com/Chocbanana/Research-Becker-Group}{online repository}} \\
	Department of Applied Mathematics\\
    University of Colorado Boulder\\
    Boulder, CO, 80309\\
	\texttt{bhavana.jonnalagadda@colorado.edu} \\
	\And
	\href{https://orcid.org/0000-0002-1932-8159}{\includegraphics[scale=0.06]{logo/orcid.pdf}\hspace{1mm}Stephen Becker} \\
	Department of Applied Mathematics\\
    University of Colorado Boulder\\
    Boulder, CO, 80309\\
	\texttt{stephen.becker@colorado.edu} \\
}
\else
\usepackage{authblk}

\newbox{\orcid}\sbox{\orcid}{\includegraphics[scale=0.06]{logo/orcid.pdf}} 
\author[1]{%
	\href{https://orcid.org/0009-0009-1054-0995}{\usebox{\orcid}\hspace{1mm}Bhavana Jonnalagadda\thanks{\texttt{bhavana.jonnalagadda@colorado.edu}}}%
}
\author[1]{%
	\href{https://orcid.org/0000-0002-1932-8159}{\usebox{\orcid}\hspace{1mm}Stehpen Becker\thanks{\texttt{stephen.becker@colorado.edu}}}%
}
\affil[1]{Department of Applied Mathematics, University of Colorado Boulder, Boulder, CO 80309}
\fi

\renewcommand{\headeright}{Research Article Preprint}
\renewcommand{\undertitle}{}
\renewcommand{\shorttitle}{\textit{arXiv}}

\hypersetup{
pdftitle={Evaluation of data driven low-rank matrix factorization for accelerated solutions of the Vlasov equation},
pdfauthor={Bhavana Jonnalagadda, Stephen Becker},
}

\begin{document}
\maketitle

\begin{abstract}
Low-rank methods have shown success in accelerating simulations of a collisionless plasma described by the Vlasov equation, but still rely on computationally costly linear algebra every time step. We propose a data-driven factorization method using artificial neural networks, specifically with convolutional layer architecture, that trains on existing simulation data. At inference time, the model outputs a low-rank decomposition of the distribution field of the charged particles, and we demonstrate that this step is faster than the standard linear algebra technique. Numerical experiments show that the method effectively interpolates time-series data, generalizing to unseen test data in a manner beyond just memorizing training data; patterns in factorization also inherently followed the same numerical trend as those within algebraic methods (e.g., truncated singular-value decomposition). However, when training on the first 70\% of a time-series data and testing on the remaining 30\%, the method fails to meaningfully extrapolate. Despite this limiting result, the technique may have benefits for simulations in a statistical steady-state or otherwise showing temporal stability.
\end{abstract}


\section*{Introduction}

\subsection*{Vlasov simulation}

Plasma describes a variety of physical phenomenon, ranging from stars to human-created fusion as part of inertial and magnetic  confinement fusion (ICF and MCF, resp.) efforts. Due to the range of forces and space- and time-scales involved, numerical simulation is challenging, and consequently there are a range of physical models that trade off accuracy for computation speed. In this work, we focus on the Vlasov equation which approximates the plasma via a \emph{density function} and ignores collisions. Instead of working with discrete particles, the density function is now a solution to a partial differential equation (PDE), which is ultimately discretized in phase space (space and velocity) and time.

The Vlasov-Maxwell model is as follows:
\begin{equation}
    \frac{\partial f}{\partial t}  + {\bf v} \cdot \nabla_{\bf x} f + \frac{q}{m}({\bf E} + {\bf v} \times {\bf B}) \cdot \nabla_{\bf v} f = 0,
    \label{eq:vlasov-complete}
\end{equation}
\begin{align}
    &-\frac{1}{c^2} \frac{\partial  \bf E}{\partial t} +\nabla \times \bf B = \mu_0  \bf J,\quad  \frac{\partial  \bf B}{\partial t} +\nabla \times \bf E =  0,
    \label{eq:maxwell}\\
    &\nabla \cdot \bf E =\frac{\rho}{\varepsilon_0},  \quad\nabla \cdot \bf B = 0,\notag
\end{align}
where $c$ is the speed of light, $\varepsilon_0$ and $\mu_0$ are the vacuum permittivity and permeability, respectively, $m$ is the mass of the particle species, $q$ is the charge of the species, and $\mathbf{E}$ and $\mathbf{B}$ are the electric and magnetic fields \cite{guoLowRankTensor2022}. 
The function to solve for is $f(t, {\bf x}, {\bf v})$, the probability distribution of the charged species, i.e., the electrons.
As the Vlasov equation is a six-dimensional equation in phase space (3 dimensions of space and 3 dimensions of velocity) plus one dimension in time, the computational cost associated with the dimensionality is a key obstacles for realistic simulations. 

In the limit of the Vlasov-Maxwell equation as the magnetic field goes to zero, ignoring relativistic effects, and using standard tricks (such as writing the electric field in terms of potentials and exploiting vector calculus identities), we arrive at the Vlasov-Poisson set of equations:
\begin{equation}
    \frac{\partial f}{\partial t}  +  {\bf{v}} \cdot \nabla_{\bf{x}}  f + {\bf{E}} (t, {\bf{x}}) \cdot \nabla_{\bf{v}}  f = 0,
    \label{eq:vlasov1}
\end{equation}
\begin{equation}
     {\bf E}(t, {\bf x}) = - \nabla_{\bf x} \phi(t, {\bf x}),  \quad -\triangle_{\bf x} \phi (t, {\bf x}) = {{\bf \rho} (t, {\bf x})} - \rho_0,
    \label{eq:poisson}
\end{equation}
where $\phi$ is the self-consistent electrostatic potential. 

On a discretized spatial and velocity grid, the numerical representation of $f(t, {\bf x}, {\bf v})$ is a 6D tensor. Recent work on low-rank representations in plasma, \cite{guoLowRankTensor2022}, has exploited the fact that if this tensor has low-rank (in a Tucker, CP or hierarchical format), then applying finite difference operators is computationally cheap.  We elaborate below on the details of this in the special case of 1 dimensions of position and 1 dimension of velocity (``1D1V''), and that is also the setting of our numerical experiments, but we note that our proposed data-driven method would actually be the most useful in the full 3D3V case since the alternative multi-linear algebra techniques for tensors are not as well-established as standard (matrix) linear algebra.

Specializing to the 1D1V case (cf.~Eq. (2.1) and (2.2) in \cite{guoLowRankTensor2022}), at a given time $t$, the numerical solution cam be represented by a matrix, with $m$ rows representing the discrete spatial locations and $n$ columns representing discrete velocity locations.  Representing this matrix as $U \in \mathbb{R}^{m\times n}$, applying the spatial derivative operator $\partial f/ \partial x$ (represented discretely as $\mathcal{D}_x\in\mathbb{R}^{m\times m}$) is equivalent to the matrix multiplication
\begin{equation}\label{eq:Dx}
    \mathcal{D}_x X
\end{equation}
and similarly applying the velocity derivative is 
\[
\mathcal{D}_v X^\top = ( X \mathcal{D}_v^\top )^\top.
\]

If $X$ has a rank-$r$ factorization $X = UV^\top$ with $U\in\mathbb{R}^{m\times r}$ and $V\in\mathbb{R}^{r\times n}$
then we can calculate 
\begin{equation}\label{eq:DxUV}
\mathcal{D}_x X = (\mathcal{D}_x U)V.
\end{equation}
Because of sparsity in the finite difference stencil, the cost of applying $\mathcal{D}_x$ to a single column of $X$ is $\mathcal{O}(m)$, and thus the cost of computing \eqref{eq:Dx} is $\mathcal{O}(mn)$. On the other hand, the cost of computing $\mathcal{D}_x U$ from \eqref{eq:DxUV} is $\mathcal{O}(mr)$, a significant savings. Likewise, the velocity derivative can be computed in $\mathcal{O}(nr)$ rather than $\mathcal{O}(mn)$. The paper \cite{guoLowRankTensor2022} gives the rest of the details for how to efficiently use these factored derivative terms. They also give empirical results showing that for factorizations based on the SVD, the approximation error can be made small enough (by choosing large enough $r$) to have good results while still leading to a computational speedup. Despite this speedup, the SVD step is still costly, and motivates our current work for finding a faster alternative.

\subsection*{Finding low-rank factorizations}

Low-rank matrix factorization (LRMF) is a powerful tool for approximating high-dimensional data by representing it as the product of two or more smaller matrices. This decomposition not only reduces the dimensionality of the data but also helps to uncover latent structures as in principle component analysis (PCA), enabling applications such as data compression, noise reduction, and feature extraction. LRMF is particularly useful in fields where the underlying data exhibits low-rank properties, such as natural language processing, computer vision, and scientific computing; by reducing data dimensionality, low-rank representations can accelerate evaluations in applications such as spectral learning and numerical solutions to dynamical systems, where computations scale with data size \cite{janzaminSpectralLearningMatrices2019}. In this particular case, approximation through LRMF affords a computational speedup to otherwise prohibitively expensive solving of the Vlasov equation.

The Singular Value Decomposition (SVD) is one of the most widely used methods for performing LRMF, offering provably optimal solutions, especially when considering a low-rank or truncated rank solution \cite{martinssonRandomizedNumericalLinear2021}. In particular, it gives the best low-rank approximation with respect to any unitarily invariant norm, such as the spectral or Frobenius norm.

However, despite its strengths, the SVD is computationally expensive for large-scale datasets, prompting research into approximate or alternative methods that could match or surpass the SVD’s performance while offering faster execution. In the case of Eq.~\eqref{eq:vlasov1}, the $x,v$ matrix is potentially exceedingly large in the respective dimensions; and when considering the full 3D3V model of Eq.~\eqref{eq:vlasov-complete}, the resulting used matrices are six-dimensional. Approximating the SVD through machine learning models like neural networks is a particularly intriguing avenue, as such methods have the potential to learn efficient low-rank representations without explicitly computing the singular values and vectors.

\subsection*{Deep learning related work}

Recent advancements in matrix factorization have leveraged deep learning to expand the scope of traditional LRMF techniques. For example, algorithms such as deep robust PCA use neural networks to solve LRMF problems for specific types of matrices such as positive semi-definite matrices. One specific example is \cite{herreraDeniseDeepRobust2023}, which claims a promising speedup although they require constraints on the input matrix structure and only present experiments for $20\times 20$ sized matrices. More broadly, most deep matrix factorization methods impose constraints on the output matrices, such as non-negativity or sparsity, to tailor the factorization for specific applications \cite{dehandschutterSurveyDeepMatrix2021, dziugaiteNeuralNetworkMatrix2015, yangOrthogonalNonnegativeMatrix2021}.

Deep matrix factorization has also been inspired by the hierarchical feature extraction capabilities of deep learning. This approach involves extracting multiple layers of features, as seen in models like DANMF \cite{danmf}, which alternates between layers of the factorized matrices for tasks such as multi-view clustering \cite{zhaoMultiViewClusteringDeep2017}. Applications of these techniques extend to hyperspectral un-mixing, recommender systems, community detection, and time-series feature extraction, demonstrating their versatility \cite{dehandschutterSurveyDeepMatrix2021, senThinkGloballyAct2019, liuTimeSeriesAnalysisLowRank2019}.

However, most of these methods focus on feature extraction or specific constraints, rather than general-purpose numerical approximations for use in dynamical equation solving. To date no work has explored deep learning for unconstrained LRMF, where both the input and output matrices are free from imposed properties, for approximating data used in solving physical dynamical equations.  Furthermore, these data-driven LRMF rely on the idea of having a dataset of similar matrices. No data-driven LRMF method could hope to rival the SVD for an arbitrary matrix.

This paper aims to fill this lack of work on the unconstrained case, presenting an unconstrained neural network-based approach to LRMF that demonstrates its utility in the context of plasma simulations. The time-series nature of these simulations gives a natural notion of a dataset from which we can train, and then predict low-rank factorizations on unseen matrices that come from the same time series.

\subsection*{Outcome reporting bias}

We address an issue stemming from the use of deep learning, particularly in the development of novel methods: outcome reporting bias, which often skews the perception of progress. Though acknowledged as an issue in many fields \cite{munafo2017manifesto}, the quick and widespread adaption of deep learning to tenuously adjacent areas of science has accelerated this problem. This bias manifests in several ways, such as selective reporting of favorable results, reliance on weak or inadequately tuned baselines, or overemphasis on niche scenarios where a proposed method performs well. These practices can lead to an overoptimistic narrative about the general effectiveness of new techniques, obscuring their limitations and making direct comparisons across studies difficult.

Weak baselines are a particularly pressing issue; for instance, it is common to see neural network-based methods evaluated without careful tuning of the baseline methods or without acknowledging that certain tasks inherently favor traditional approaches due to their simplicity or inherent mathematical properties. This creates a distorted picture where new methods appear more advantageous than they truly are when applied in broader or more realistic contexts. Within the more specific context of deep learning for solving of fluid-related partial differential equations, it has been found that outcome reporting bias and publication bias are widespread \cite{mcgreivy2024weak}.

In order to resist this problem, we aim to not participate in selective reporting, spin bias, or other flavors of outcome reporting bias. We instead present a balanced evaluation of the proposed method used within this paper, acknowledging and displaying the limitations found; additionally, we use a rigorous baseline against which the proposed method is evaluated in comprehensive ways.

\subsection*{Contributions}

This paper introduces ConvMF, a convolutional neural network-based approach to low-rank matrix factorization (LRMF), tailored to handle structured time-series data generated by high-dimensional plasma simulations governed by the Vlasov equation. By leveraging the representational power of neural networks, ConvMF provides a novel alternative to traditional methods such as SVD, demonstrating its ability to capture key features of complex datasets while achieving competitive performance in reconstruction accuracy and (after amortizing the training phase) outperforming SVD in execution time. Thorough analysis of the loss and output from the model provides insight into the use of deep learning for interpolation, extrapolation, and its applicability to LRMF for both stable and unstable temporal timeseries data. Additionally, within the analysis we aim to present honest appraisal, contributing to the ongoing discourse about reproducibility and rigor in machine learning research. By explicitly discussing the limitations and avoiding overoptimistic claims, this paper seeks to provide a foundation for future work that is grounded in realistic expectations and a clear understanding of where neural network-based methods can and cannot excel in comparison to traditional approaches, in the context of LRMF for accelerated simulation solving.

\section*{Materials and methods}

Precise details on implementations and usage of statistical methods, built models, data generation code, and evaluation metrics are located at \nameref{S1_Code}.

\subsection*{Data}

The data used in this study consisted of simulated numerical solutions to the 1D1V Vlasov-Poisson equation, using the code written by Guo and Qiu \cite{guoLowRankTensor2022}. For a given simulation (specified by a standard initial condition), the dataset is a timeseries of $m\times n$ matrices, where $m$ is the number of spatial points and $n$ is the number of velocity points.

\begin{figure}[!h]
\begin{adjustwidth}{-2.25in}{0in}
\centering
\includegraphics[width=.3\linewidth]{"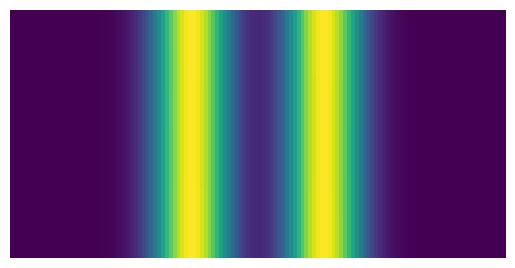"}\hfill
\includegraphics[width=.3\linewidth]{"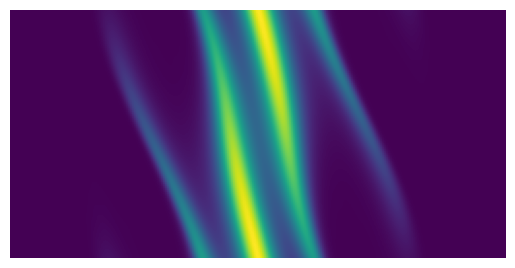"}\hfill
\includegraphics[width=.3\linewidth]{"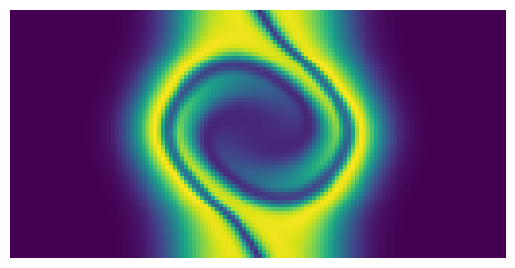"}
\caption{{\bf Visualization of generated data.}
The phase-space visual representation of the generated plasma data used, at size $128\times256$ pixels, with selected frames (matrices) from the full timeseries sequence. The $x$-axis is space and $y$-axis is velocity. {\bf Left}: One of the first frames in the sequence, and also the ``easiest'' frame for the networks to learn (lowest average loss over all networks). {\bf Middle}: A randomly selected frame from the middle of the sequence. {\bf Right}: One of the last frames, and also the frame that was ``hardest" to learn (highest average loss over all networks).}
\label{fig:dataimgs}
\end{adjustwidth}
\end{figure}

Two matrix sizes were analyzed: $64\times128$ and $128\times256$. Larger sizes, such as $256\times512$, as well as 4D versions from 2D2V simulations, were excluded due to being prohibitively expensive given computational constraints. The initial condition for the $64\times128$ matrices corresponded to a strong 1D problem (order 5), while for the $128\times256$ matrices, the initial condition was a two-stream instability 1D problem (order 5). Representative data for the $64\times128$ case is shown in Fig~\ref{fig:dataimgs}, where the $x$-axis is the spatial dimension and the $y$-axis is the velocity dimension.

Most of our experiments focused on interpolating, hence we \emph{randomly} partitioned the timeseries dataset into 70\% training data and 30\% validation/test data. Validation data was used during training to monitor and inform on loss, while the test data remained unseen until the final evaluation. We also looked at extrapolation, performed by training on the first 70\% of the time-series and testing/validating on the remaining 30\%.

Additionally, matrices generated with randomized initial conditions (e.g., smoothed noise without constraints), using the same code as above, were included in an attempt to improve model generalization during extrapolation tests.

\subsection*{Neural network structure and experimental configurations}

The neural networks employed in this study were designed to perform low-rank matrix factorization by outputting two matrices, $U$ and $V$, whose product approximates the input matrix; see Fig.~\ref{fig:nnstructure} for the network architecture that was ultimately selected. Separate networks were trained for the different input sizes. The networks were built with the following design principles:

\begin{itemize}
    \item \textbf{Architecture}: The networks featured an initial sequence of linear layers (referred to as the ``stem'') that branched into two parallel paths of linear layers (``fork'') to generate the $U$ and $V$ outputs.
    \item \textbf{Variations}: Some network configurations incorporated convolutional layers, while others relied solely on fully connected layers. Various number of layers was used. Prior work \cite{dehandschutterSurveyDeepMatrix2021} had found that three hidden layers was typically sufficient for matrix factorization tasks; our tests never looked beyond 6 layers given this obvservation.
    \item \textbf{Time-invariance}: Recurrent layers and sequence-input (such as that for transformer attention layers) were deliberately excluded to focus on time-invariant decomposition without leveraging prior temporal information.
    \item \textbf{Output ranks}: Networks were trained separately for various output ranks, without employing transfer learning, to allow comparative evaluation of consistent architectures.
    \item \textbf{Hyperparameters}: Multiple configurations of activation functions, optimizers, learning rates, and layer compositions were tested.
    \item \textbf{Ablation studies}: Network components were systematically removed or modified to assess their contributions.
    \item \textbf{Comparison model}: A re-implementation of the DAMNF model \cite{danmf} was included for comparative purposes, the results of which are discussed in \nameref{S1_Appendix}.
\end{itemize}

\subsection*{Evaluation and comparison}

The primary evaluation metric was the normalized reconstruction loss, which is commonly used by low rank matrix factorizations done by neural networks \cite{dehandschutterSurveyDeepMatrix2021}, and calculated as follows:
\begin{equation} 
\label{eq:loss}
\min_{\substack{U \in \mathbb{R}^{m\times r} \\ V \in \mathbb{R}^{r\times n}}} \dfrac{\|X - UV\|^2_F}{\|X\|^2_F} ,
\end{equation}
where $\|\cdot\|_F$ is the Frobenius norm, i.e., $\|X\|_F^2 = \sum_{i=1}^m\sum_{j=1}^n X_{ij}^2$.
The scaling by $1/\|X\|_F^2$ was only used during
performance evaluation but not during training.

Additional evaluation was used to assess performance.
In terms of both accuracy and execution time, comparisons of the proposed method were made with the standard linear algebra technique (the SVD). The SVDs were implemented by calling the standard \texttt{scipy.linalg} and \texttt{scipy.sparse.linalg} libraries in Python.
Output validity was determined by finding the loss from using a given (fixed) output from the model (either $U$ or $V$) and then using least-squares to find the computed reciprocal output. Additionally, the execution time of the trained models on test data was measured.

\section*{Results and discussion}

\begin{figure}[!h]
\centering
\includegraphics[width=1.0\textwidth]{"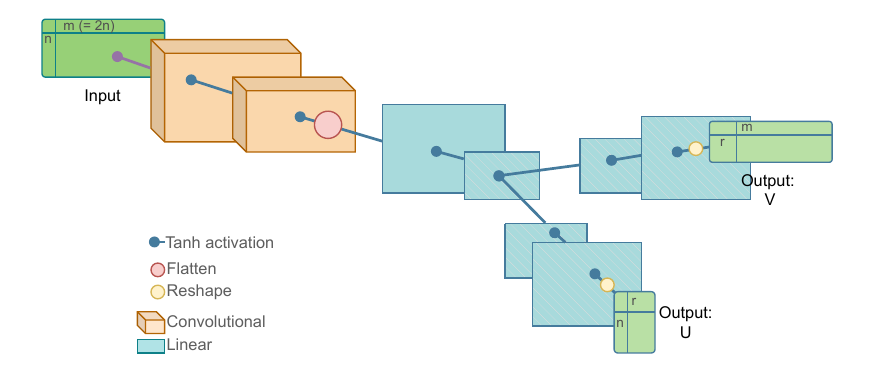"}
\caption{{\bf ConvMF structure.}
The structure of the selected best network after training and hyperparameter testing. The input goes through several convolutional and then linear layers, before splitting into two paths (fork) with several linear layers for each output, U and V. Visual sizes of the layers are not to actual scale of the dimensions used, and are instead to indicate the relative scale of the layers to each other.}
\label{fig:nnstructure}
\end{figure}

Various architectures and hyperparameters were tested, with the options tested listed in Table~\ref{tab:hyperparameters}. The tested architectures were also tested with various output ranks ($5$--$30$), and with different input matrix sizes. The selected architecture and model, hereby ``ConvMF'', was chosen from having the lowest loss across all tested ranks and input sizes, on validation data. The architecture of the chosen best model can be seen in Fig.~\ref{fig:nnstructure}. The network consists of 2 convolutional layers, feeding into 2 linear layers (all with activation function between), at which point there is a ``fork'' of 3 linear layers each for the two outputs $U$ and $V$. The training and validation loss curve of ConvMF can be seen in \nameref{S1_Fig}.

\begin{table}[h]
    \centering
    \caption{{\bf Tested hyperparameters.} The following options were tested in various combinations, and the final used hyperparameters in ConvMF are bolded. All were also potentially tested with output of various ranks ($5$--$30$), and with various input matrix sizes. $n,m=$ input matrix dimensions $n\times m$, $r=$ rank of output U and V.}
    \begin{tabular}{ll}
    \toprule 
    Hyperparameter & Values tested \\
    \midrule 
    \textbf{Layer activation} & \{ReLU, \textbf{Tanh}, LeakyRelu, Sigmoid\} \\
    \textbf{Optimizer} & \{\textbf{ADAM}, SGD, Adagrad\} \\
    \textbf{Learning rate} & \{0.001, 0.0005, \textbf{0.0001}\} \\
    \textbf{Number of conv layers} & \{1, \textbf{2}, \ldots, 5\} \\
    \textbf{Conv parameters}  &  Kernel Size=\{\textbf{3}, \textbf{5}, 6\}, Stride=\{\textbf{1}, 2, 3\}, \\
        & Padding=\{\textbf{0}, 1, 2, \textbf{3}\}, Dilation=\textbf{1} \\
    \textbf{Number of stem layers} & \{1, \textbf{2},\ldots, 6\} \\
    \textbf{Stem layer dim} & \big\{100, \textbf{200}, \textbf{500}, 1000, 2000, $\dfrac{n\cdot m}{2^{\{1,\ldots,5\}}}$\big\} \\
    \textbf{Number of fork layers} & \{1, \textbf{2},\ldots, 6\} \\
    \textbf{Fork layer dim} & \big\{100, \textbf{200}, \textbf{300}, \ldots, 1000, $\dfrac{m\cdot r}{9-2^{\{1,\ldots,5\}}}$\big\} \\
    \bottomrule
    \end{tabular}
    \label{tab:hyperparameters}
\end{table}

During training, it was found that networks with more layers and of larger dimensional size (8+ layers) performed comparatively and competitively against those with less and smaller layers, with final validation loss being within $10^{-4}$ of each other; thus, we chose the smaller network with same other hyperparameters. Additionally, networks with the convolutional layers only benefited from about 2 layers, with more (and larger convolutional layer output) actually increasing loss, likely due to the increased difficulty in training them.

\subsection*{Accuracy on interpolation data}
\begin{figure}[h]
\centering
\includegraphics[width=1.0\textwidth]{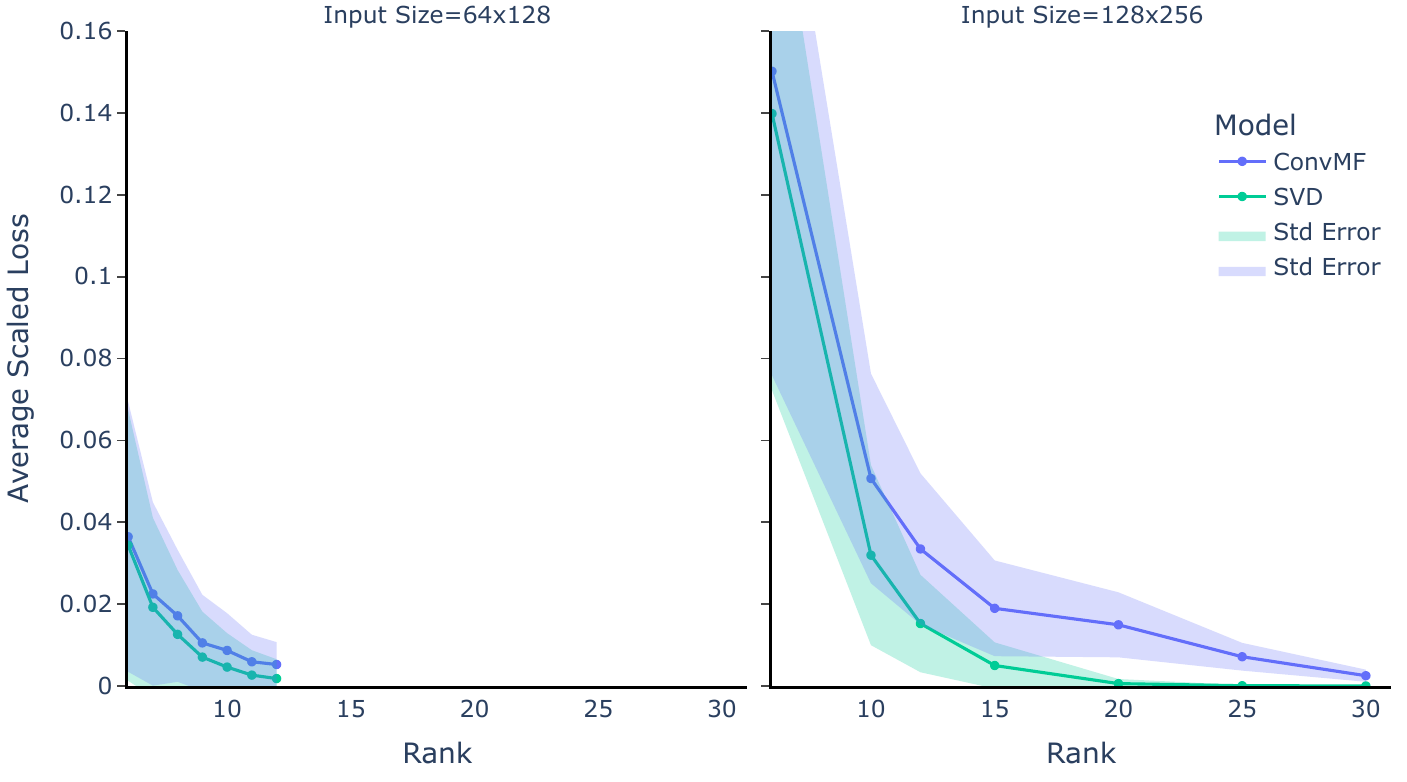}\hfill
\includegraphics[width=1.0\textwidth]{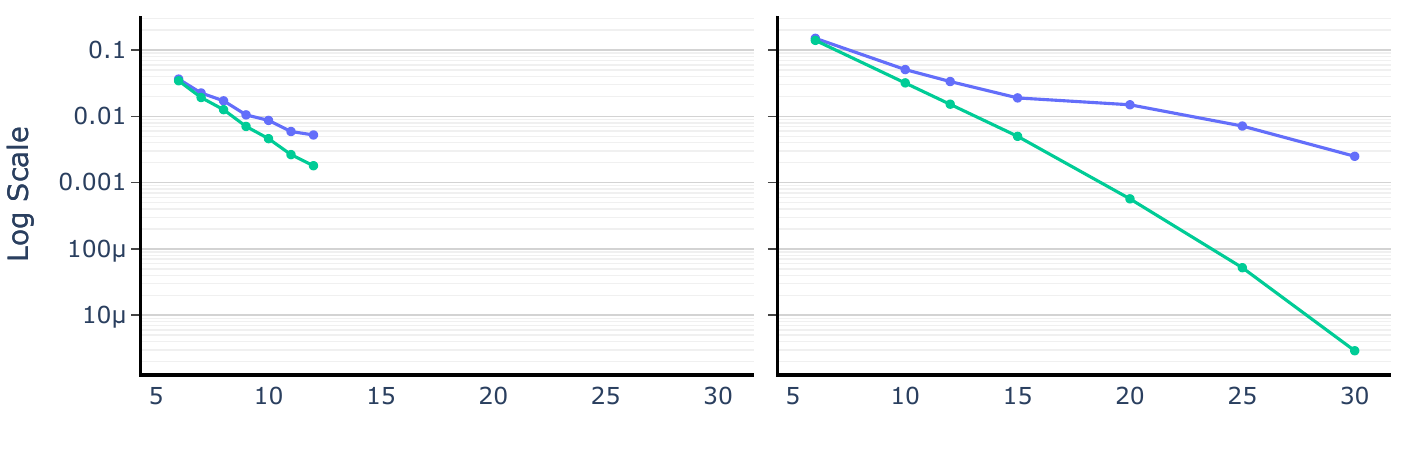}
\caption{{\bf Rank vs average scaled loss.}
The average scaled loss, on validation and held-out test data, for the best neural network tested (``ConvMF'') and for the SVD, across different ranks of the resulting $U$ and $V$. Different ranks for ConvMF were done with separately trained networks for each rank.
{\bf Left}: The resulting loss per rank for input data of size $64\times128$. 
{\bf Right}: The resulting loss per rank for input data of size $128\times256$. 
{\bf Bottom}: The same two plots in log scale for $y$.}
\label{fig:rankavgloss}
\end{figure}

The rank of the output, versus the average scaled loss, on validation and test data, can be seen in Fig~\ref{fig:rankavgloss}, showing both a linear and a log scale. Loss was scaled by the Frobenius norm as discussed in Methods. 
While the SVD has smaller loss, as expected due to its optimality, the ConvMF loss is within comparable range.
Across a variety of ranks, as a rule-of-thumb, the error of the SVD at rank $r$ can be achieved by ConvMF using rank $\approx r+2$.
For very high accuracies (i.e., very high rank), the SVD starts to significantly outperform ConvMF, but at medium and low accuracy the difference is smaller.

\begin{figure}[h!]
\begin{adjustwidth}{-2.25in}{0in}
\includegraphics[width=\linewidth]{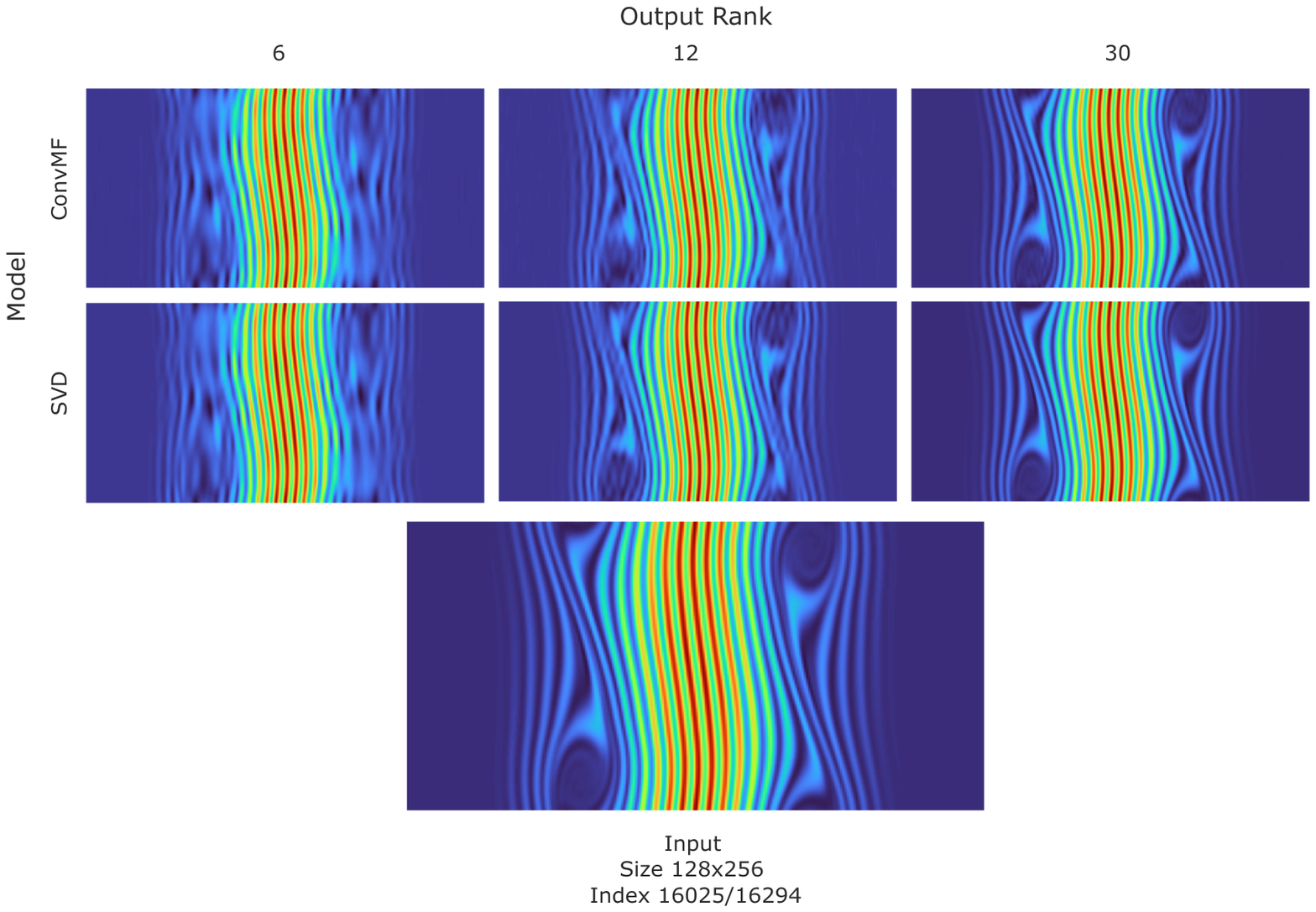}
\caption{{\bf Reconstructed output from models vs original input.}
The display of matrices formed from reconstruction (matrix multiplication) from the constituent output, shown against the original input of size $128\times256$, for selected ranks. The specific input at time index 16025 was chosen as it had the highest average loss across all models (and ranks) used. {\bf Top row}: The reconstruction from the output $U, V$ from ConvMF, at ranks 6, 12, 30. {\bf Middle row}: The reconstruction from the output $U, \Sigma, V$ from SVD, at ranks 6, 12, 30. {\bf Bottom image}: The original input to both methods, at size $128\times256$ and 98\% into the timeseries.}
\label{fig:outputimgs}
\end{adjustwidth}
\end{figure}

In Fig.~\ref{fig:outputimgs}, we present a visual comparison between the reconstructed outputs of ConvMF and SVD for the \(128 \times 256\) input that exhibited the highest average loss across all tested models and ranks. The figure shows outputs at ranks 6, 12, and 30, providing a progression of how reconstruction quality improves as the rank increases. Both methods appear visually indistinguishable at corresponding ranks, highlighting the effectiveness of ConvMF in capturing the core features of the dataset similarly to SVD.

For the $128\times256$ matrix at time index 16025 shown in the figure, 
at rank 30, which corresponds to \(23\%\) of the original input size, the reconstructed outputs for both ConvMF and SVD are nearly indistinguishable from the original matrix, even upon close inspection. This emphasizes the capacity of these low-rank methods to preserve fine-grained details despite the substantial dimensionality reduction. The ability of ConvMF to match SVD in reconstructing such complex data reinforces its potential utility in handling high-dimensional data with intricate structures.

While we do not propose low-rank factorization as a replacement for more advanced image compression techniques, such as those optimized for natural image datasets, this visualization demonstrates the ability of ConvMF to efficiently capture and represent key features within the data. It also underscores the interpretability and adaptability of ConvMF for scientific applications where preserving critical features at reduced dimensions is essential.

\subsection*{Execution time and complexity}

An advantage of ConvMF is that after training, using the network to provide a low-rank factorization is computationally cheap. To quantify this, we compare against the time to compute a low-rank factorization via the SVD. We use two standard SVD solvers, the first being the default SVD from Python's \texttt{scipy.linalg.svd} which uses standard \texttt{LAPACK} dense linear algebra routines. For a $m\times n$ matrix with $m\le n$, these routines cost $\mathcal{O}(m^2n)$ flops. In the plots, this algorithm is referred to as ``Basic SVD''.

The second SVD solver is the Krylov based solver in \texttt{scipy.sparse.linalg.svds} which uses standard \texttt{ARPACK} libraries. Each iteration requires a matrix-vector multiplication, with a cost of $\mathcal{O}(mn)$ flops. The number of required multiplications depends on properties of the matrix and the desired tolerance, but it is typically $\mathcal{O}(r)$, and hence the overall complexity is $\mathcal{O}(mnr)$ which is faster than the dense SVD when $r \ll \min\{m,n\}$ \cite{scipy, arpack}. In the plots, this algorithm is referred to as ``Faster SVD''.

\begin{figure}[h]
\centering
\includegraphics[width=0.8\textwidth]{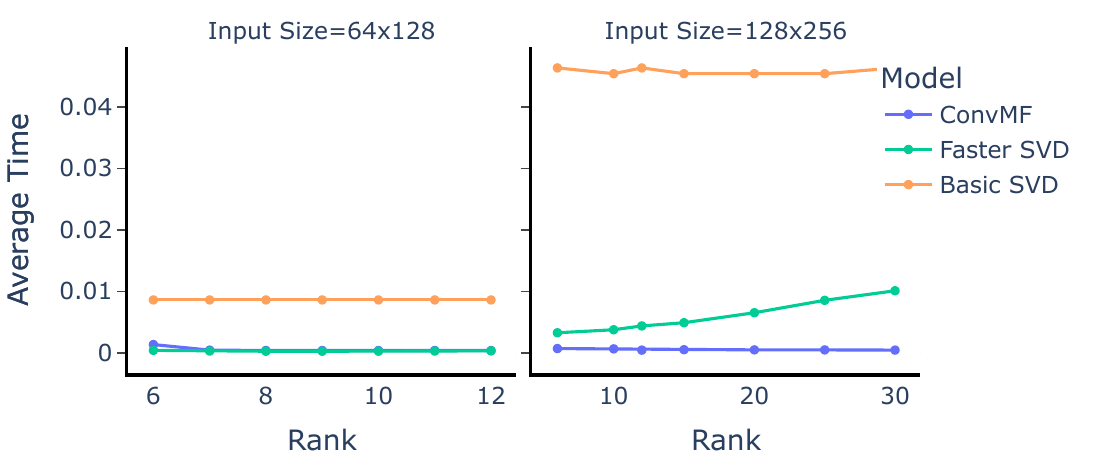}
\caption{{\bf Average execution time across ranks.}
The averaged execution time for the two SVD methods tested, along with the selected network, evaluating on test data. {\bf Left:} The times for input matrix size $64\times128$. {\bf Right:} The times for input matrix size $128\times256$.}
\label{fig:exectime}
\end{figure}

For the both input sizes, Fig.~\ref{fig:exectime} shows the average execution time across ranks. For the smaller input, execution time against the Faster SVD had negligible difference; for the larger input, ConvMF is relatively faster, with the speed difference becoming greater as rank increased since computational cost of the Faster SVD increases linearly with rank, as expected from our above analysis.

To predict a general model of the runtime cost of ConvMF, 
we fix the factors that are invariant and static in its structure across different ranks and input sizes; namely, the convolution stride, filter/kernel size, number of filters, number of channels, and the dense/linear layers' dimensions (excluding the last layer output dimension). Given this, and the fact that the last linear layer's execution time scales linearly in relation to the output rank, we can express the dominating factors in the time complexity of ConvMF as $\mathcal{O}(C + mn + r)$ \cite{dumoulinGuideConvolutionArithmetic2018}, where C is the constant factor affecting the time the whole network takes to evaluate (and where the effect of C is greater than that of $mn$ at scales used in this paper). Output rank scales the time separately from the input size, and thus contributes to a near-flat line of execution time for ConvMF in Fig.~\ref{fig:exectime}.

The above model of the runtime for ConvMF at larger sizes and ranks is predicated on maintaining good accuracy, which is not something we can guarantee \emph{a priori}, but the experiments we do have of doubling the input size and using a variety of ranks are promising and suggest it may be possible to have rapid runtime. The caveat is that ConvMF needs a comparatively long training time, but it is conceivable that a base model could be trained and then modifications made via transfer learning, thus greatly reducing offline training cost.

\subsection*{Loss analysis}

\begin{figure}[H]
\centering
\includegraphics[width=\textwidth]{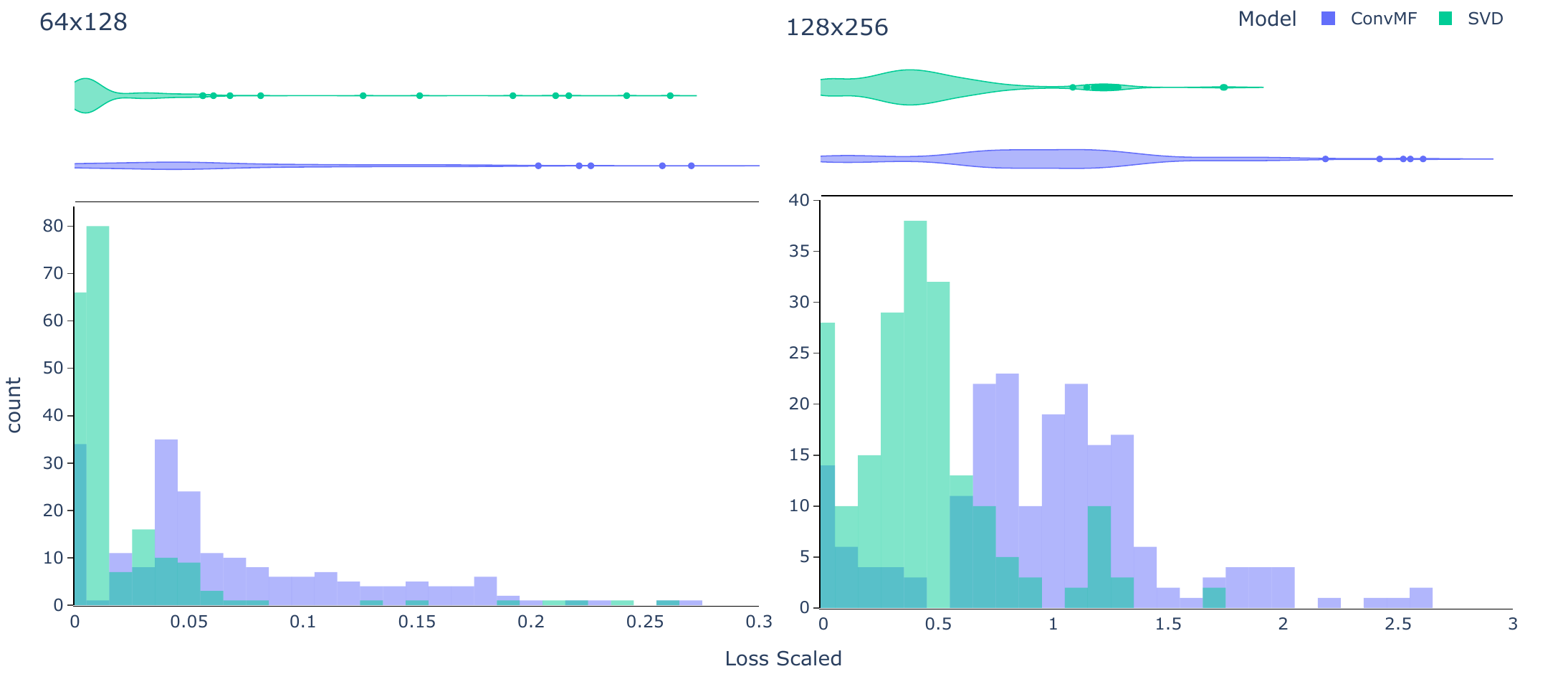}
\caption{{\bf Loss by input matrix for rank 12.}
Histograms of the loss per each input matrix in the timeseries test + validation data, for output rank 12. {\bf Left}: Input matrix size $64\times128$. {\bf Right}: Input matrix size $128\times256$.}
\label{fig:losshist}
\end{figure}

The histogram of loss, per each input matrix from the time series, for rank 12 and for two input sizes, can be seen in Fig.~\ref{fig:losshist}. In general, the spread of loss from ConvMF was greater than that of the SVD methods, though also following the the general distribution of the latter.

\begin{figure}[H]
\centering
\includegraphics[width=0.9\textwidth]{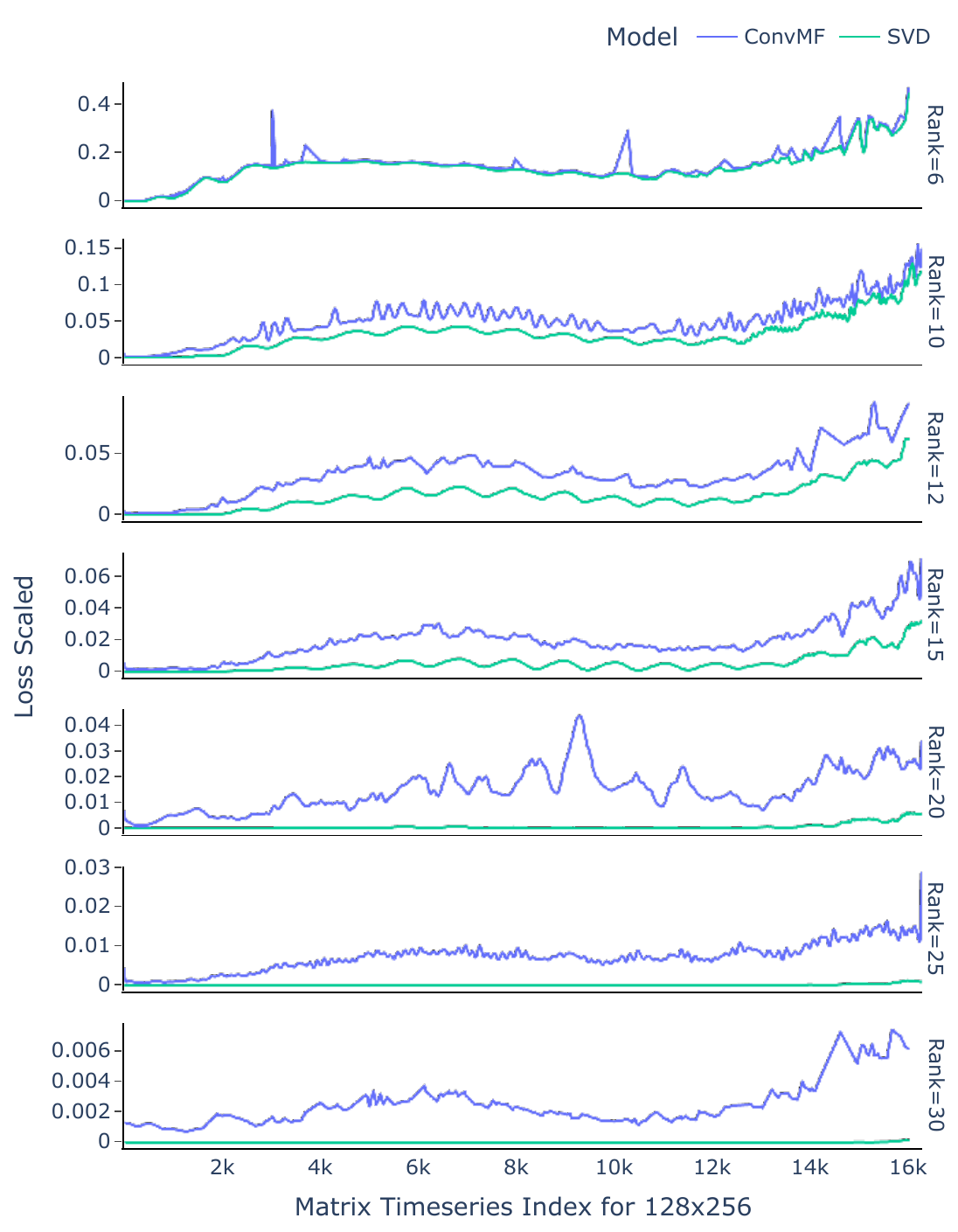}
\caption{{\bf Scaled loss vs input matrix of size $128\times256$.}
The scaled loss, for each input matrix in the test and validation data, ordered by the matrix timeseries index, for input size $128\times256$, shown for each tested rank.}
\label{fig:lossbyinput}
\end{figure}

This effect can also be seen,  with more illuminating detail, in Fig.~\ref{fig:lossbyinput} where the loss per each input in the timeseries is shown, for various ranks and for input size $128\times256$, all on unseen validation and testing data. We can observe that there is an overall trend for increase in the loss for both the neural network and SVD method as the time series progresses, with a bump in the loss around 1/3 into the timeseries; the trends also closely follow each other, though this is less obvious at higher ranks, where the loss is low enough that the SVD line trends flat. This suggests that there are inherent characteristics in the data that complicate the factorization process at ranks low enough to have reduced included information, and that the difficulty of the of the problem of LRMF has a hard bottom with this data. 

Notably, a distinct rise in loss is observed approximately one-third of the way through the time series, after which the loss continues to escalate. This pattern suggests that the data becomes progressively more challenging to factorize as it evolves, indicating the presence of complex underlying structures that neither method fully captures at lower ranks. The close alignment of the loss trends between ConvMF and SVD implies that both models encounter similar difficulties in representing the data's latent structures, particularly when rank constraints limit the amount of information retained. This parallel behavior underscores the inherent characteristics of the dataset, which likely contains temporal or structural dependencies that complicate factorization, especially when reduced-rank approximations are used. Such challenges could be attributed to increasing heterogeneity or noise within the data, which both methods struggle to model effectively under constrained settings.

\subsection*{Error analysis}

ConvMF produces an approximation $U\cdot V$ of the original matrix $X$. In this section, we look at what part of the approximation dominates the  error. Letting $\text{col}(X)$ denote the column space of a matrix $X$ (the span of the columns), we can ask if the error is primarily due to (1) a mismatch between $\text{col}(X)$ and $\text{col}(U)$, or (2) a mismatch between $\text{col}(X^\top)$ and $\text{col}(V^\top)$, or (3) correct row and column spaces but incorrect coefficients inside these spaces.

To examine each of the above 3 causes, given the outputs $U,V$ from ConvMF, we calculate three additional losses that postprocess some or all of the ConvMF outputs with additional linear algebra techniques:

\begin{enumerate}
    \item $\min_{ \widetilde{U}\in \mathbb{R}^{m \times r} } \frac{ \| X - \widetilde{U}V \|_F }{\|X\|_F}$, ``Calculated $U$''
    \item $\min_{ \widetilde{V}\in \mathbb{R}^{r \times n} } \frac{ \| X - U\widetilde{V} \|_F }{\|X\|_F}$, ``Calculated $V$''
    \item $\min_{ \widetilde{\Sigma}\in \mathbb{R}^{r \times r} } \frac{ \| X - U\widetilde{\Sigma}V \|_F }{\|X\|_F}$, ``Calculated $\Sigma$''.
\end{enumerate}

All three of the above problems can be rewritten as least-squares problems and solved via standard software.

To give an example of drawing inference from the metrics, if, say, the loss from ``Calculated $U$'' matched the loss from the SVD, then we can conclude that most of the error in ConvMF is due to having the wrong $U$, and in particular the wrong column space.

\begin{figure}[H]
\begin{adjustwidth}{-2.25in}{0in}
\begin{subfigure}{0.5\linewidth}
\caption{}
\includegraphics[width=0.9\linewidth]{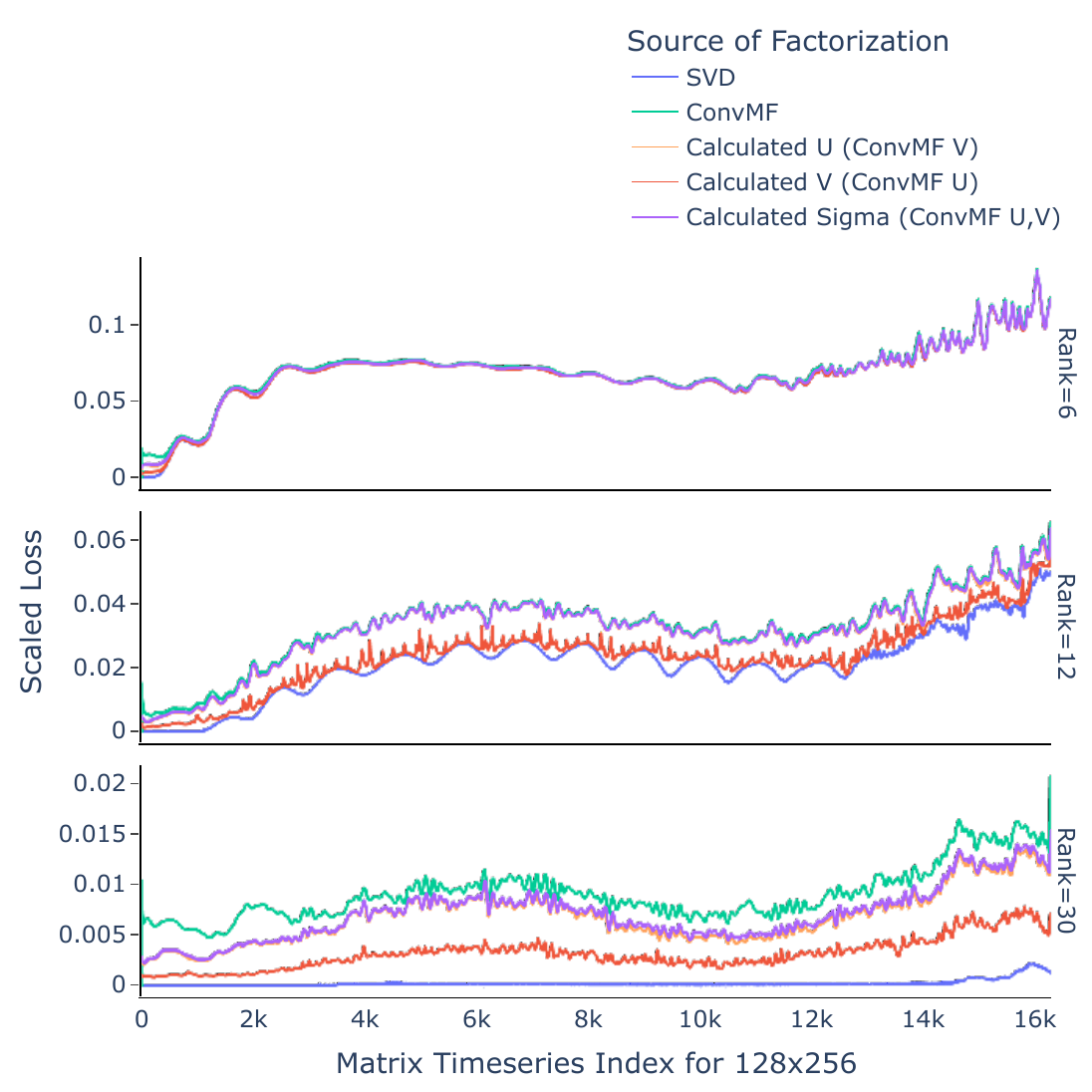}
\end{subfigure}
\begin{subfigure}{0.5\linewidth}
\caption{}
\includegraphics[width=0.9\linewidth]{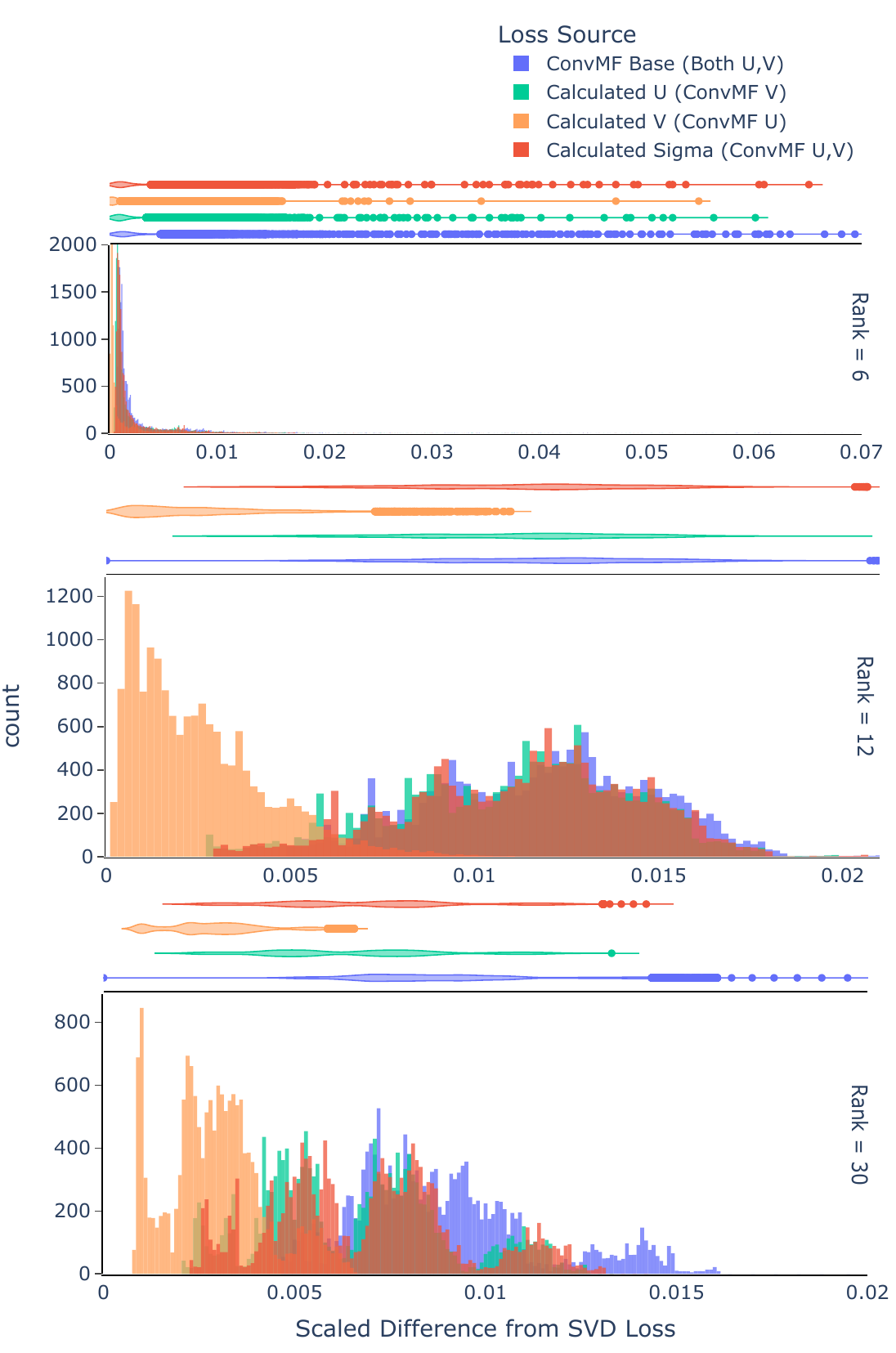}
\end{subfigure}
\caption{{\bf Loss difference between input and calculated sources, by output rank.}
The scaled loss difference between a calculated source (or the multiplied output from models ConvMF, SVD), and the input. The calculated \{U, V, S\} indicates that the reciprocal output from ConvMF is fixed, and then used in a least-squares calculation (using the defined loss function and given input) to find the optimal calculated value. {\bf A}: The calculated scaled loss difference by timeseries index of the input, by rank. {\bf B}: Histogram of the difference between the loss from SVD, and the loss from using the calculated value, by rank.}
\label{fig:correctness}
\end{adjustwidth}
\end{figure}

The results, illustrated in Fig.~\ref{fig:correctness}, highlight distinct differences between the ConvMF outputs. Notably, the calculated $V$ matrix — obtained by fixing the ConvMF output $U$ and solving for $V$ — produced losses closest to those of SVD, outperforming other configurations by a significant margin. 
Some of this difference compared to solving for $U$ may simply be because $n=2m$ so $V$ has twice as many degrees of freedom as $U$. However, looking at Part B of Fig.~\ref{fig:correctness}, the difference between $U$ and $V$ appears to be more than just a factor of $2$.

Additionally, the calculated $\Sigma$, derived by solving while fixing both $U$ and $V$, displayed remarkable similarity to the calculated $U$ matrix. This observation suggests that the $U$ output from ConvMF closely resembles the corresponding output from SVD, to a much greater extent than $V$.

The discrepancy between the quality of $U$ and $V$ outputs is not fully explained but may be attributed to the inherent differences in their structures and optimization. The $U$ matrix is half the size of $V$, potentially making it easier for the network to optimize; additionally, the backward propagation of parameter updates, driven by loss minimization, may affect the weights corresponding to the fork output of $U$ more significantly than those for $V$. Or, the discrepancy could be because position and velocity are fundamentally different quantities.

Overall, these findings demonstrate that ConvMF can achieve a remarkably close approximation to SVD, particularly for lower ranks, where the differences between the two methods become comparable. This also suggests that a hybrid data-drive method plus linear-algebra post-processing method is an interesting future topic to explore.

\subsection*{Extrapolation testing}

\begin{figure}[H]
\begin{adjustwidth}{-2.25in}{0in}
\begin{subfigure}{0.5\linewidth}
    \begin{subfigure}{1.0\linewidth}
    \caption{}
    \includegraphics[width=1.0\linewidth]{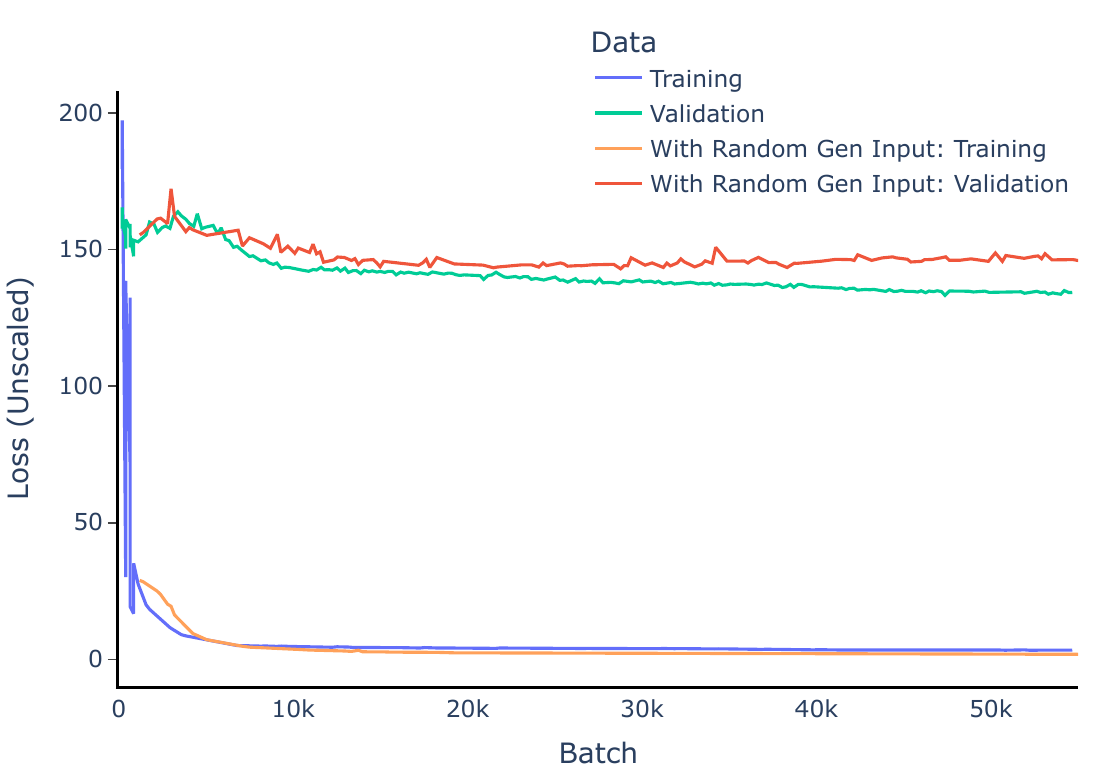}
    \label{fig:tvloss-seq}
    \end{subfigure}
    \begin{subfigure}{1.0\linewidth}
    \caption{}
    \includegraphics[width=1.0\linewidth]{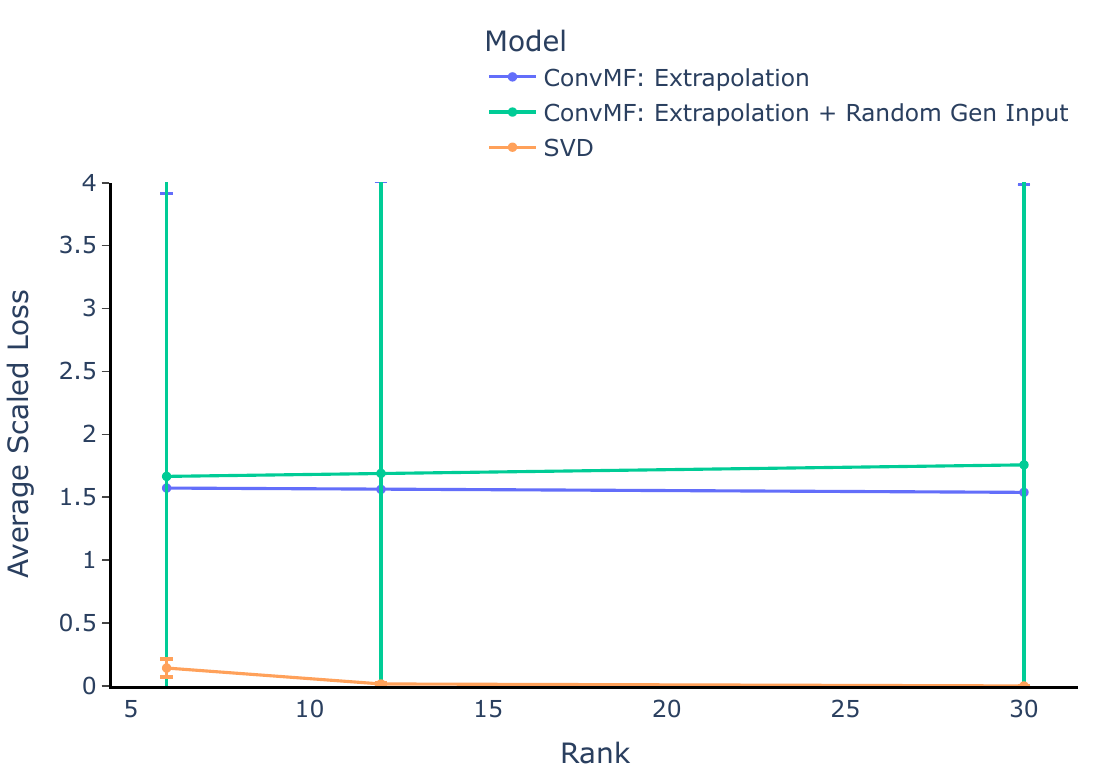}
    \label{fig:rankavgloss-seq}
    \end{subfigure}
\end{subfigure}
\begin{subfigure}{0.5\linewidth}
    \caption{}
    \includegraphics[width=1.0\linewidth]{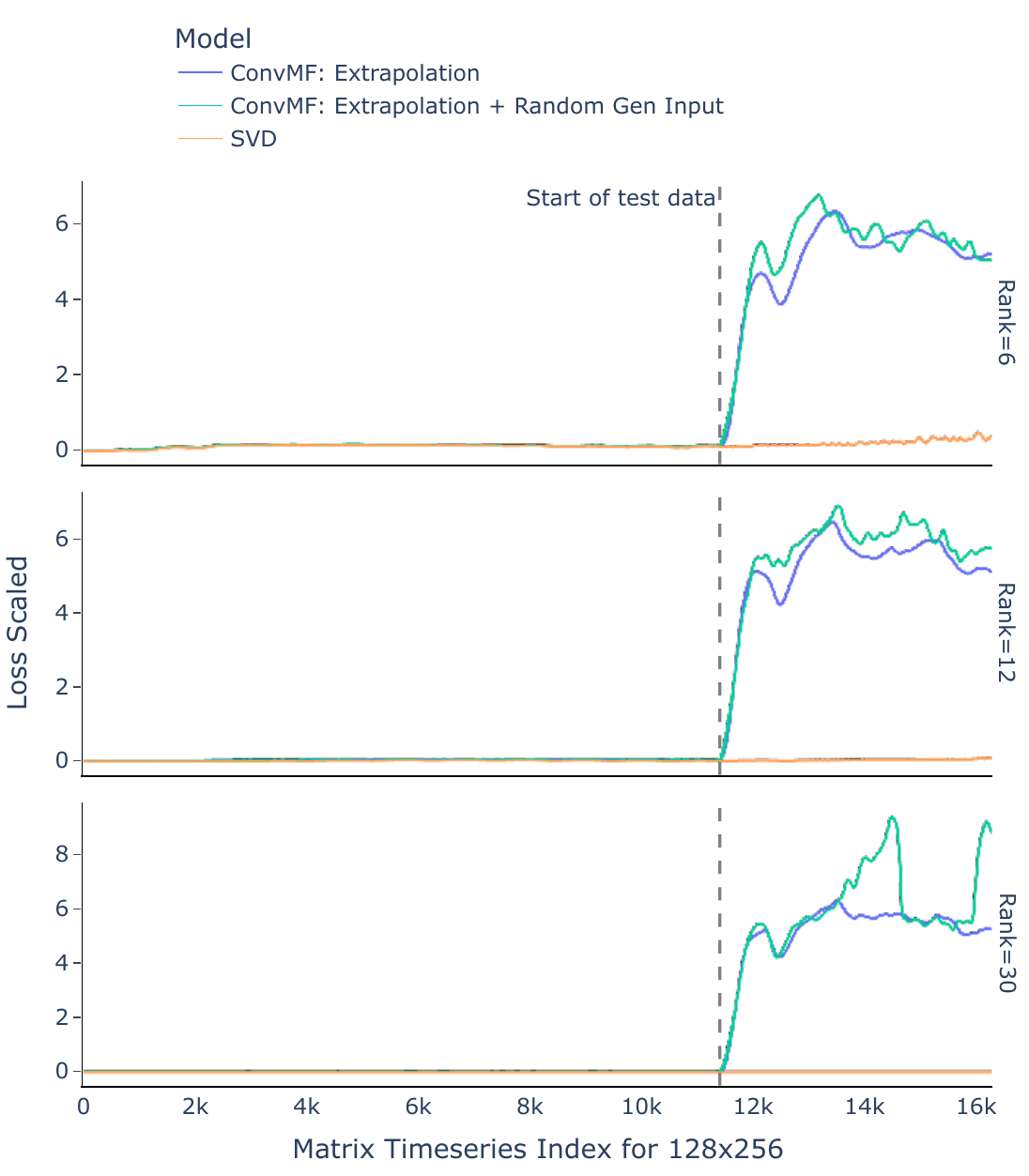}
    \label{fig:lossbyinput-seq}
\end{subfigure}
\caption{{\bf Extrapolation testing results.} Results for the best network tested (ConvMF), where it was trained on the first 70\% of data and the validation/test data was the last 30\% of the data. The inclusion of randomly generated data was for training on only, and validation was still performed only on the 30\% of original plasma simulation data. {\bf A}: The training and validation loss curve for input size $128\times256$ and output rank 12.  {\bf B}: The average scaled loss, by rank. {\bf C}: The scaled loss, for each input matrix in the test and validation data, ordered by the matrix timeseries index, for input size $128\times256$, separated by rank. The dashed line indicates the end of data seen by the models during training.}
\label{fig:extrap}
\end{adjustwidth}
\end{figure}

Testing was done for extrapolation, where the network was trained on the first 70\% of data (randomly fed) and testing and validation was performed on the last 30\% of data. Additionally, extrapolation testing was performed by including randomly generated data (from different initial conditions) in training, created as discussed in Methods. The training and validation loss curves can be seen in Fig.~\ref{fig:extrap}.A; the validation curves never lowered below a large threshold above the training curves, with the inclusion of randomly generated data slightly worsening the loss. This outcome can also be seen in Fig~\ref{fig:extrap}.B, with the average loss staying consistently at a high threshold above that achieved by SVD. Similarly, in Fig~\ref{fig:extrap}.C, we observe that once the training data ends and the test data begins, ConvMF struggles with maintaining effective factorization, and the addition of randomly generated data appears to slightly exacerbate this challenge.

Several insights into the behavior of ConvMF when faced with unseen data can be drawn. Importantly, ConvMF's consistent ability to capture temporal structures within the training data demonstrates its potential for learning stable representations under interpolation conditions, which could form the basis for improvements in extrapolation. Additionally, the introduction of randomly generated data, while increasing loss, indicates that ConvMF is sensitive to input variability and can adapt its factorization process to diverse data patterns. This adaptability suggests that further fine-tuning of network architecture or training procedures might enable ConvMF to leverage such random data more effectively, potentially improving its generalization to future sequences. These results highlight opportunities for refining training techniques to better capture temporal evolution, which is critical for applications requiring robust extrapolation capabilities. The inclusion of randomly generated data during training not improving the model's ability to generalize, and in some cases worsening the performance, suggests that ConvMF may be over-fitting to the training data and failing to capture some of the underlying patterns in the later part of the time series. Further work will be necessary to address this limitation, potentially through regularization techniques or alternative training strategies.

\section*{Conclusion}

\subsection*{Future directions}

Future work will aim to address the limitations observed in this study and expand the applicability of ConvMF for broader use cases. A promising direction includes higher-dimensional data, by applying ConvMF to 4D or 6D tensor datasets; though previously constrained by memory limitations, it could provide insights into its performance on more complex structures. Furthermore, the classical baseline for tensor data is not as well understood theoretically and is also slower to compute numerically, so there is more room for improvement over standard methods. For the ConvMF model, there is no obstacle, other than increases in training time and memory, for increasing the dimensionality.

Incorporating temporal dependencies, such as through recurrent or attention-based layers, may improve the model’s ability to handle extrapolation and generalize across diverse datasets; though negating the generalization abilities or testing for such, furthermore increasing execution time, it might be beneficial in use cases different from the one posed in this paper. Leveraging domain knowledge through physics-informed neural networks (PINNs) could enhance ConvMF’s performance on physics-based datasets, such as those governed by the Vlasov equation or similar models; utilizing the specialized regularization techniques in PINNs could better fit ConvMF to the need and inherent features present within a given dataset.

\subsection*{Summary}

This study explored the use of deep learning (ConvMF) for low-rank matrix factorization of time-series plasma data. Through extensive experimentation, we demonstrated that ConvMF can effectively give low-rank factors for reconstructing complex matrices, at a competitive accuracy and reconstruction to the SVD, particularly at lower ranks. The results highlight ConvMF's ability to capture key features of the data and maintain competitive performance in terms of reconstruction loss, showing promise as a method for learning low-rank matrix factorizations, especially when applied to datasets that exhibit temporal stability.

The method was shown to adapt well to structured time-series data, effectively replicating the results of SVD under interpolation conditions. However, when trained with additional randomly generated data, ConvMF struggled to generalize, often leading to worsened performance. This highlights the model's tendency to learn specific patterns inherent to the dataset rather than generalizable features. The inability of ConvMF to generalize well in extrapolation testing underscores the need for better handling of temporal dependencies, suggesting avenues for architectural improvements. Additionally, the observation that ConvMF approximates SVD closely at lower ranks highlights its potential for specific, constrained use cases rather than as a general replacement for the SVD.

Importantly, ConvMF performed competitively with SVD in both accuracy and execution time for larger datasets for interpolation, demonstrating potential scalability advantages. When leveraged for long sequence datasets, with large input size and temporally stable values, deep learning shows the ability to be efficient in time and for factorizing with interpolation. 

In summary, ConvMF offers a computationally efficient approach to LRMF for temporally stable data and lays the groundwork for future innovations in matrix factorization techniques tailored to scientific datasets. By addressing current limitations and pursuing the outlined directions, ConvMF can evolve into a more robust and versatile tool for high-dimensional efficient solving of plasma simulations.

\clearpage

\section*{Acknowledgments}

We acknowledge and thank Jing-Mei Qiu for extending the information and code needed to generate the simulated plasma data used within this paper.
Both authors gratefully acknowledge support from the DOE under award DE-SC0023346: Center for Hierarchical and Robust Modeling of Non-Equilibrium Transport (CHaRMNET).

\nocite{*}

\clearpage

\section*{Supporting information}

\paragraph*{S1 Fig.}
\label{S1_Fig}
{\bf Training/Validation loss curve.} The training and validation loss curve for the best network tested (ConvMF), with loss being unscaled by norm.

\begin{figure}[!h]
\centering
\includegraphics[width=1.0\textwidth]{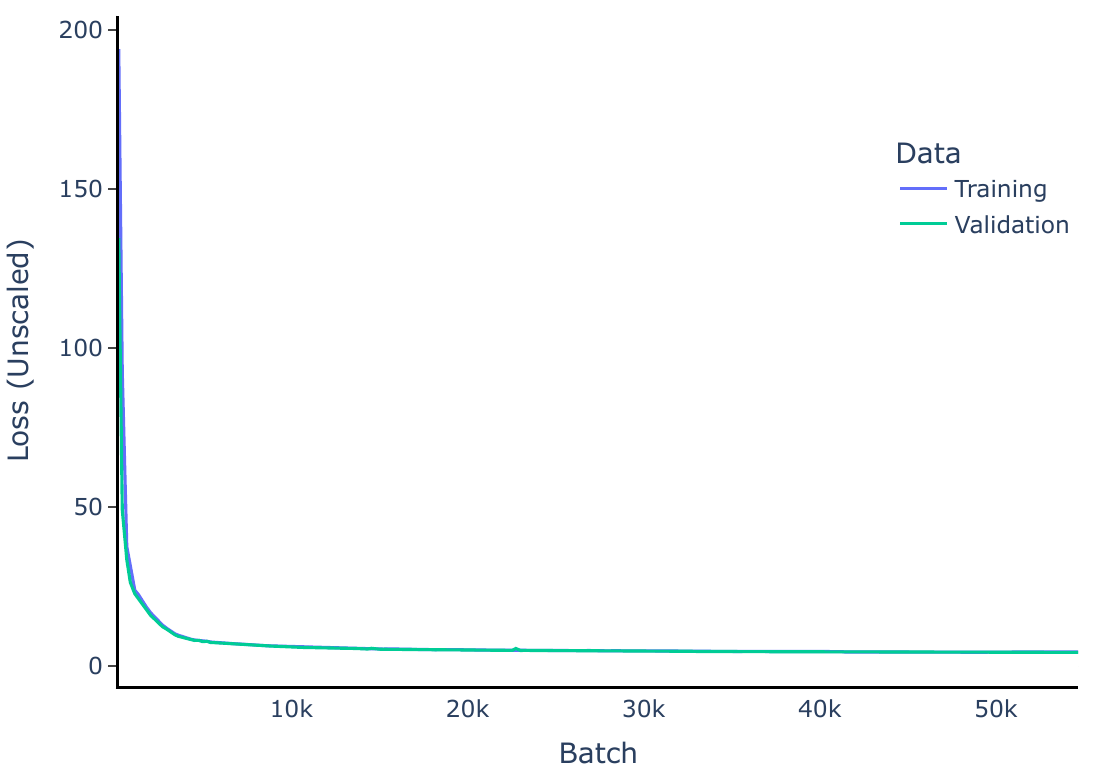}
\end{figure}

\paragraph*{S1 Appendix.}
\label{S1_Appendix}
{\bf Comparison again re-implemented method.} DANMF \cite{danmf}, a network structure of alternating layers for matrix factorization where "factorizations of the basis and coding matrices are alternated along layers", was re-implemented in the chosen framework. The loss was, at best, an order of $10^2 - 10^4$ times worse than the average performance achieved by models implemented in this paper, for both input sizes. This potentially suggests an issue in the re-implementation or application to the simulated plasma data used.

\paragraph*{S1 Code Link.}
\label{S1_Code}
{\bf Source code for research performed.} The code used to generate the data, some of the generated data itself, and all model-related code (training, inference, analysis, comparison) is publicly available on Github: \url{https://github.com/Chocbanana/Research-Becker-Group}.

\end{document}